\documentclass[12pt]{article}
\usepackage[ansinew]{inputenc}
\usepackage{comment}
\usepackage{amsmath}
\usepackage{tikz}
\usepackage{amssymb}
\usepackage{amsthm}
\usepackage{fancyvrb}
\usepackage{authblk}
\usepackage{hyperref}
\usepackage{ifpdf}

\newcounter{Theorem}
\newenvironment{Theorem}[1][]{\refstepcounter{Theorem}\par\medskip
   \noindent \textbf{Theorem~\theTheorem. #1} \rmfamily}{\medskip}
   \newcounter{Lemma}
\newenvironment{Lemma}[1][]{\refstepcounter{Lemma}\par\medskip
   \noindent \textbf{Lemma~\theLemma. #1} \rmfamily}{\medskip}
\newcounter{Definition}
\newenvironment{Definition}[1][]{\refstepcounter{Definition}\par\medskip
   \noindent \textbf{Definition~\theDefinition. #1} \rmfamily}{\medskip}
\newcounter{Example}

\newcounter{Remark}

\newcounter{Corollary}
\newenvironment{Corollary}[1][]{\refstepcounter{Corollary}\par\medskip
   \noindent \textbf{Corollary~\theCorollary. #1} \rmfamily}{\medskip}

\title{Notes on the Neighborhood Polynomials}
   
\makeatletter
\renewcommand\@date{{%
  \vspace{-\baselineskip}%
  \large\centering
  \begin{tabular}{@{}c@{}}
    Maryam Alipour\\
    \normalsize malipour@hs-mittweida.de
  
  \end{tabular}

  \bigskip

  \textsuperscript{}University of Applied Sciences Mittweida

  \bigskip

  \date{}
}}
\makeatother

\author{}

\begin{document}
\maketitle

	\begin{abstract}
	The \emph{neighborhood polynomial} of graph $G$, denoted by $N(G,x)$, is the generating function for 
		the number of vertex subsets of $G$ which are subsets of open neighborhoods of vertices in $G$. For any graph polynomial, it can be useful to generate a new family of polynomials by introducing some restrictions and characterizations. In this paper, we investigate two new graph polynomials that are obtained from $N(G,x)$ by adding \textit{independence} or \textit{connectivity} restrictions to the vertex subsets or to the subgraphs induced by the vertex subsets which are generated by $N(G,x)$. These new polynomials are not only related to $N(G,x)$, but also having strong connections to other known graph polynomials of $G$ or its subgraphs, such as \textit{independence polynomials} or \textit{subgraph component polynomials}.
	\end{abstract}
\noindent
\textbf{Keywords:} \hspace{2mm} Neighborhood Polynomials, Domination Polynomials, Independence Polynomials, Subgraph Component Polynomials.

\noindent
\textbf{2010 Mathematics Subject Classification:} \hspace{2mm} 05C69, 05C31.

\section{Introduction}
\noindent
Let $G=(V,E)$ be a simple, finite and undirected graph. For a vertex $v\in V$, the \emph{open neighborhood} of $v$, denoted by $N(v)$, is defined by $N(v)=\{u\mid \{u,v\}\in E\}$.

\vspace{2mm}
\noindent
The \emph{neighborhood complex} of $G$, denoted by $\mathcal{N}(G)$ and introduced in \cite{Lov78}, is the family of all subsets of open neighborhoods of 
vertices of $G$,
$$ 
\mathcal{N}(G) = \{A \mid A\subseteq V , \exists v\in V : A\subseteq N(v) \}.
$$

\noindent
The \emph{neighborhood polynomial} of graph $G$, denoted by $N(G,x)$ and introduced in \cite{BN08}, is the ordinary generating function for the 
neighborhood complex of $G$.
\begin{equation*}
N(G,x) = \sum_{U\in \mathcal{N}(G)} x^{|U|}.
\label{eq:1}
\end{equation*}

\noindent
The neighborhood polynomial of graph $G$ can be represented as $N(G,x)=\sum_{k=0}^{n} n_{k}(G)x^k$ where $n$ is the order of $G$ and $n_{k}(G)$ is the number of vertex subsets of cardinality $k$ in $G$ which are subsets of open neighborhoods of vertices in $G$. In other words $n_k(G)=|\{A \mid A\in \mathcal{N}(G),|A|=k\}|$.

\vspace{2mm}
\noindent
A \emph{dominating set} of a graph $G$ is a vertex subset $W$ of $G$ such that the closed neighborhood of $W$ is equal to V, where the \emph{closed neighborhood} of $W\subseteq V$ is defined by $N[W] = \bigcup_{w\in W}N(w) \cup W$. If the family of all dominating sets of $G$ is denoted by $\mathcal{D}(G)$, the \emph{domination polynomial} of a graph, introduced in \cite{AL00}, can be defined as the ordinary generating function for $\mathcal{D}(G)$. In other words,
\[
D(G,x) = \sum_{W\in\mathcal{D}(G)}x^{|W|}.
\]
As it is proved in \cite{HT17}, the following relation between the domination polynomial of $G$ and the neighborhood polynomial of the complement of $G$, $N(\bar{G},x)$, is satisfied which motivates the study of neighborhood polynomials.
$$D(G,x) + N(\bar {G},x) = (1+x)^{|V|}.$$

\vspace{3mm}

\noindent
In Section \ref{sect:ind} of this paper, we investigate how the addition of the \emph{independence} restriction to the vertex subsets in $\mathcal{N}(G)$ changes the neighborhood polynomial of $G$. We give the title of \emph{independent neighborhood polynomial} to this new graph polynomial and then explore some of its properties, the effect of several graph operations on it and study its relation with other known graph polynomials.

\vspace{2mm}
\noindent
In Section \ref{sect:con}, we consider the addition of \emph{connectivity} restriction to the subgraphs induced by the vertex subsets in $\mathcal{N}(G)$ and call the new polynomial \emph{connected neighborhood polynomial}. Then, we investigate the properties of this polynomial and its characteristics as well.

\section{Independent Neighborhood Polynomial} \label{sect:ind}
\noindent
The \textit{independence polynomial} of graph $G$, denoted by $I(G,x)$ and introduced in \cite{GH83}, is the polynomial
\begin{equation*}\label{def:independence polynomial}
I(G,x)=\sum_{k=0}^{n}i_{k}(G)x^{k},
\end{equation*}
where $n$ is the order of $G$ and $i_{k}(G)$ is the number of independent vertex subsets of cardinality $k$ in $G$.

\vspace{2mm}
\noindent
Here, we introduce the \textit{independent neighborhood polynomial} of $G$ which is the ordinary generating function for the number of subsets of open neighborhoods of vertices of $G$ which form independent sets in $G$.

\vspace{2mm}
\begin{Definition}
The \textit{independent neighborhood polynomial} of graph $G$ is the polynomial
\begin{equation*}\label{def:independent neighborhood polynomial}
N^{(i)}(G,x)=\sum_{k=0}^{n}n_{k}^{(i)}(G) x^{k},
\end{equation*}
where $n$ is the order of $G$ and $n_{k}^{(i)}(G)$ is the number of independent vertex subsets of cardinality $k$ which are subsets of open neighborhoods in $G$. 
\end{Definition}

\vspace{2mm}
\noindent
If a graph $G$ is a tree, then every vertex set with a common neighbor in $G$ also forms an independent set in $G$, otherwise $G$ contains a cycle as a subgraph which contradicts the fact that $G$ is a tree. This argument leads to the following result.
\begin{Lemma}
If $G$ is a tree, then
\begin{equation*}
    N^{(i)}(G,x)=N(G,x).
\end{equation*}
\end{Lemma}

\noindent
Also, if $G$ is a cycle of order $n$ with $n>3$, clearly it does not contain any cycle of order $k$, with $k<n$, as a subgraph which means any vertex subset with a common neighbor in $G$ forms an independent set in $G$. This argument leads to the following result.
\begin{Lemma}
If $G$ is a cycle of order $n$ with $n>3$, then
\begin{equation*}
    N^{(i)}(G,x)=N(G,x).
\end{equation*}
\end{Lemma}

\noindent
The \textit{disjoint union} of two graphs $G_{1}=(V_{1}, E_{1})$ and $G_{2}=(V_{2}, E_{2})$ with disjoint vertex sets $V_{1}$ and $V_{2}$, denoted by $G_{1}\cup G_{2}$ is a graph with the vertex set $V_{1}\cup V_{2}$, and the edge set $E_{1}\cup E_{2}$.

\vspace{2mm}
\begin{Theorem}
Suppose $G_{1}$ and $G_{2}$ are vertex disjoint graphs. The independent neighborhood polynomial of the disjoint union $G_{1}\cup G_{2}$ of these two graphs satisfies
     \begin{equation*}
         N^{(i)}(G_{1}\cup G_{2},x)=N^{(i)}(G_{1},x)+N^{(i)}(G_{2},x)-1.
     \end{equation*}
\end{Theorem}
\begin{proof}
It is clear that every independent vertex subset with a common neighbor in $G_{1}$ forms also an independent vertex subset with a common neighbor in $G_{1}\cup G_{2}$. The same argument holds for $G_{2}$. The only set that is counted more than once is the empty set and this over-counting is corrected by subtracting $1$.
\end{proof}

\noindent
The \textit{join} of two graphs $G_{1}=(V_{1}, E_{1})$ and $G_{2}=(V_{2}, E_{2})$ with disjoint vertex sets $V_{1}$ and $V_{2}$, denoted by $G_{1}+G_{2}$, is the disjoint union of $G_{1}$ and $G_{2}$ together with all those edges that join vertices in $V_{1}$ to vertices in $V_{2}$.
\vspace{2mm}
\begin{Theorem}
Suppose $G_{1}=(V_{1}, E_{1})$ and $G_{2}=(V_{2}, E_{2})$ are vertex disjoint graphs. The independent neighborhood polynomial of the join $G_{1}+G_{2}$ of these two graphs satisfies
\begin{equation*}
N^{(i)}(G_{1}+G_{2},x)=I(G_{1},x)+I(G_{2},x)-1.
\end{equation*}
\end{Theorem}
\begin{proof}
Let $X$ be an independent vertex subset of $G_{1}+G_{2}$ with a common neighbor in $G_{1}+G_{2}$. We can distinguish the following cases:
\begin{enumerate}
    \item[] - Suppose $X\subseteq V_{1}$ or $X\subseteq V_{2}$. Since in $G_{1}+G_{2}$ there is an edge between each vertex in $V_{1}$ and every vertex in $V_{2}$, the vertices in every independent subset of $V_{1}$ have a common neighbor in $G_{2}$. Similarly, the vertices in every independent subset of $V_{2}$ have a common neighbor in $G_{1}$. All such subsets are generated by $I(G_{1},x)+I(G_{2},x)-1$ in which $I(G,x)$ is the independence polynomial of graph $G$. By subtracting $1$, we correct the over-counting of the empty set.
    
    \item[] - Since there is an edge between each vertex in $V_{1}$ and every vertex in $V_{2}$, the set $X$ can not contain vertices from both $V_{1}$ and $V_{2}$ and at the same time remains an independent set. This argument completes the proof.
\end{enumerate}
\end{proof}

\begin{Corollary}
For the given graphs $G_{1}$ and $G_{2}$ the following holds.
\begin{equation*}
N^{(i)}(G_{1}+G_{2},x)=I(G_{1}+G_{2},x).
\end{equation*}
\end{Corollary}

\noindent
The \emph{Cartesian product} of graphs $G_{1}=(V_{1},E_1)$ and $G_{2}=(V_{2},E_2)$ 
with disjoint vertex sets $V_1$ and $V_2$, denoted by $G_{1} \square G_2$, is a 
graph with vertex set $V_{1}\times V_{2}=\{(u,v) \mid u\in V_{1}, v\in V_{2}\}$ where the vertices $x=(x_{1},x_{2})$ and $y=(y_{1},y_{2})$ are adjacent in 
$G_{1} \square G_2$ if and only if $[x_{1}=y_{1}\text{ and }\{x_{2},y_{2}\}\in E_{2}]$ 
or $[x_{2}=y_{2}\text{ and } \{x_{1},y_{1}\}\in E_{1}]$.

\vspace{2mm}
\begin{Theorem}
Let $G_{1}=(V_{1},E_{1})$ and $G_{2}=(V_{2},E_{2})$ be vertex disjoint graphs with no isolated vertex. Then
\begin{align*}
 N^{(i)}(G_{1} \square G_{2},x)=& 1+|V_{1}|(N^{(i)}(G_{2},x)-1)+|V_{2}|(N^{(i)}(G_{1},x)-1)\\
 &+\sum_{(u,v)\in V_{1}\times V_{2}} (I(G_{1}[N_{G_{1}}(u)],x)-1)(I(G_{2}[N_{G_{2}}(v)],x)-1)\\
 &-|V_{1}||V_{2}|x - 2|E_{1}||E_{2}|x^{2}.
\end{align*}
\end{Theorem}
\begin{proof}
Clearly, $G_{1} \square G_{2}$ contains $|V_{1}|$ copies of $G_{2}$ and at the same time $|V_{2}|$ copies of $G_{1}$. Every independent vertex subset of each copy of $G_{1}$ with a common neighbor in that copy of $G_{1}$ can be transformed to an independent vertex subset of $G_{1} \square G_{2}$ with a common neighbor in $G_{1} \square G_{2}$ and the same argument holds for the copies of $G_{2}$. All such vertex subsets are generated by $1+|V_{1}|(N^{(i)}(G_{2},x)-1)+|V_{2}|(N^{(i)}(G_{1},x)-1)$ and since any set containing a vertex $(u,v)$ with $u\in V_{1}$ and $v\in V_{2}$ is counted twice, once in $N^{(i)}(G_{1},x)$ as the vertex $u$ in a copy of $G_{1}$ and once more in $N^{(i)}(G_{2},x)$ as the vertex $v$ in a copy of $G_{2}$, by subtracting $|V_{1}||V_{2}|x$ we correct this over-counting.

\vspace{2mm}
\noindent
By $\sum_{(u,v)\in V_{1}\times V_{2}} (I(G_{1}[N_{G_{1}}(u)],x)-1)(I(G_{2}[N_{G_{2}}(v)],x)-1)$, we generate all non-empty independent subsets of the open neighborhood of the vertex $(u,v)$ in $G_{1} \square G_{2}$ which contain at least one vertex of the form $(u^{'},v)$ with $u^{'} \in N_{G_{1}}(u)$, and at least one vertex of the form $(u,v^{'})$ with $v^{'} \in N_{G_{2}}(v)$ for all $(u,v)\in V_{1}\times V_{2}$.

\vspace{2mm}
\noindent
Every copy of an edge $\{a,b\}\in E_{1}$ together with every copy of an edge $\{x,y\}\in E_{2}$ forms a four cycle in $G_{1} \square G_{2}$ in which, each of the vertex subsets $\{(a,x),(b,y)\}$ and $\{(b,x),(a,y)\}$ are counted twice and the double-counting of such subsets is corrected by subtracting $2|E_{1}||E_{2}|x^{2}$. 
\end{proof}

\vspace{2mm}

\noindent
Let $G=(V,E)$ be a graph. For any positive integer $r$ the $r$-expansion of $G$, denoted by $exp(G,r)$, is the graph obtained from $G$ by replacing every vertex $v\in V$ with an independent set $I_{v}$ of size $r$ and replacing every edge $\{u,v\}\in E$ with a complete bipartite graph $K_{r,r}$ with the bipartite sets $I_{u}$ and $I_{v}$.
\begin{Theorem}
Let $G=(V,E)$ be a graph, $r$ be a positive integer and $exp(G,r)$ be the $r$-expansion of $G$. Then,
$$N^{(i)}(exp(G,r),x)=N^{(i)}(G,((1+x)^{r}-1)).$$
\end{Theorem}
\begin{proof}
Suppose $X$ is an independent vertex subset with a common neighbor $v$ in $G$. Then any vertex subset in $exp(G,r)$ which contains at least one vertex from at least one $I_{x}$ forms an independent set with all vertices in $I_{v}$ as common neighbors. On the other hand, any independent vertex subset with a common neighbor in $exp(G,r)$ is of this form and all such sets are generated by $N^{(i)}(G,(1+x)^{r}-1)$.
\end{proof}
\section{Connected Neighborhood Polynomial}\label{sect:con}
In this section, we introduce the \emph{connected neighborhood polynomial} of a graph which is the ordinary generating function for the number of subsets of open neighborhoods of vertices of a given graph which induce a connected subgraph.
\begin{Definition}
The \emph{connected neighborhood polynomial} of graph $G$ is the polynomial \begin{equation*}\label{def:connected neighborhood polynomial}
N^{(c)}(G,x)=\sum_{k=0}^{n}n_{k}^{(c)}(G) x^{k},
\end{equation*}
where $n$ is the order of $G$ and $n_{k}^{(c)}$ is the number of vertex subsets of cardinality $k$ which are subsets of open neighborhoods of vertices in $G$ and induce a connected subgraph in $G$.
\end{Definition}

\vspace{2mm}

\noindent
Let $G=(V(G), E(G))$ and $H=(V(H), E(H))$ be simple graphs. An \emph{isomorphism} from $G$ to $H$ is a bijection $\phi: V(G)\rightarrow V(H)$ such that $\{u,v\}\in E(G)$ if and only if $\{\phi(u),\phi(v)\}\in E(H)$. If such an isomorphism from $G$ to $H$ exists we write $G\cong H$.

\vspace{4mm}

\noindent
The following lemma can be easily proved.
\begin{Lemma}
Let $G$ be a graph of order $n$. Then the following statements are satisfied.
\begin{enumerate}
    \item[] - Suppose $K_{n}$ is the complete graph of order $n$ and $G\cong K_{n}$, then $N^{(c)}(G,x)=N(G,x)= (1+x)^{n}-x^{n}$,
    \item[] - Suppose $T_{n}$ is a tree of order $n$ and $G\cong T_{n}$, then $N^{(c)}(G,x)=1+nx$,
    \item[] - Suppose $C_{n}$ is a cycle of order $n$ with $n>3$ and $G\cong C_{n}$, then $N^{(c)}(G,x)=1+nx$.
\end{enumerate}
\end{Lemma}
\begin{Theorem}
Suppose $G_{1}$ and $G_{2}$ are vertex disjoint graphs. The connected neighborhood polynomial of the disjoint union $G_{1}\cup G_{2}$ of these two graphs satisfies
     \begin{equation*}
         N^{(c)}(G_{1}\cup G_{2},x)=N^{(c)}(G_{1},x)+N^{(c)}(G_{2},x)-1.
     \end{equation*}
\end{Theorem}
\begin{proof}
It is clear that every vertex subset with a common neighbor in $G_{1}$ which induces a connected subgraph in $G_{1}$ forms a subset with a common neighbor in $G_{1}\cup G_{2}$ which also induces a connected subgraph in $G_{1}\cup G_{2}$. The same argument holds for $G_{2}$. The only set that is counted twice is the empty set which is corrected by subtracting $1$.
\end{proof}

\vspace{4mm}
\noindent
Let $G=(V,E)$ be a simple graph of order $n$. The \emph{subgraph component polynomial} of graph $G$, denoted by $Q(G;x,y)$ and introduced in \cite{TAM11}, is defined by
\begin{equation*}
    Q(G;x,y)=\sum_{i=0}^{n}\sum_{j=0}^{n} q_{ij}(G)x^{i}y^{j},
\end{equation*}
where 
\begin{equation*}
    q_{ij}(G)=|\{X\subseteq V \hspace{1mm}|\hspace{1mm} k(G[X])=j \hspace{1mm} \wedge |X|=i  \hspace{1mm} \}|
\end{equation*}
in which $k(G[X])$ gives the number of connected components of the graph induced by the vertex subset $X\subseteq V(G)$.\\
\vspace{2mm}

\noindent
The coefficient of $y^{k}$ in $Q(G;x,y)$, written as $[y^{k}]Q(G;x,y)$ and denoted by $Q_{k}(G;x)$, is the ordinary generating function for the number of vertex subsets in $G$ that induce a subgraph in $G$ with exactly $k$ components: $Q_{k}(G;x)=[y^{k}]Q(G;x,y)$. It is clear that $Q_{1}(G;x)$ is the generating function for the number of vertex subsets that induce a connected subgraph. This particular case is renamed to \emph{subgraph polynomial} and rephrased as follows:
\begin{equation*}
    S(G,x)=Q_{1}(G;x)=\sum_{k=0}^{n}s_{k}(G)x^{k},
\end{equation*}
in which $s_{k}(G)$ gives the number of vertex subsets of size $k$ which induce a connected subgraph in $G$. For further information on subgraph polynomial, see \cite{TAM11}.
\vspace{2mm}

\begin{Theorem}
Suppose $G_{1}=(V_{1}, E_{1})$ and $G_{2}=(V_{2}, E_{2})$ are vertex disjoint graphs. The connected neighborhood polynomial of the join $G_{1}+G_{2}$ of these two graphs satisfies
\begin{align*}
    N^{(c)}(G_{1}+G_{2},x)&=S(G_{1},x)+S(G_{2},x)\\
    &+(N(G_{1},x)-1)((1+x)^{|V_{2}|}-1)\\
    &+(N(G_{2},x)-1)((1+x)^{|V_{1}|}-1)\\
    &-(N(G_{1},x)-1)(N(G_{2},x)-1).
\end{align*}
\end{Theorem}
\begin{proof}
Let $X$ be a vertex subset in $G_{1}+G_{2}$ inducing a connected subgraph in $G_{1}+G_{2}$ with a common neighbor in $G_{1}+G_{2}$. We can distinguish the following cases:
\begin{enumerate}
    \item[] - Suppose $X\subseteq V_{1}$ or $X\subseteq V_{2}$. Since any subset of $V_{1}$ that induces a connected subgraph in $G_{1}$ forms a vertex subset of $G_{1}+G_{2}$ which induces a connected subgraph in $G_{1}+G_{2}$ and obviously has all vertices in $V_{2}$ as common neighbors in $G_{1}+G_{2}$. So, we just need to count all vertex subsets in $G_{1}$ that induce connected subgraphs in $G_{1}$. The same argument holds for $G_{2}$. All such sets are generated by $S(G_{1},x)+S(G_{2},x)$.
    \item[] - Suppose $X$ contains at least one vertex from each of $V_{1}$ and $V_{2}$. Let $X$ be the disjoint union of non-empty vertex subsets $A_{1}$ and $A_{2}$ where $A_{1}\subseteq V_{1}$ and $A_{2}\subseteq V_{2}$. If there is a common neighbor $v_{1}$ for the vertices of $A_{1}$ in $G_{1}$, then after joining $G_{1}$ and $G_{2}$ all vertices in $X$ have $v_{1}$ as their common neighbor in $G_{1}+G_{2}$ and $X$ induces a connected subgraph in $G_{1}+G_{2}$. A similar arguments holds for the case that there is a vertex $v_{2} \in V_{2}$ as a common neighbor for vertices of $A_{2}$ in $G_{2}$. Such sets $X$ are generated by $(N(G_{1},x)-1)((1+x)^{|V_{2}|}-1)+(N(G_{2},x)-1)((1+x)^{|V_{1}|}-1)$. By subtracting $(N(G_{1},x)-1)(N(G_{2},x)-1)$ we correct the double-counting of those sets $X=A_{1}\cup A_{2}$ where $A_{1}\in \mathcal{N}(G_{1})$ and $A_{2}\in \mathcal{N}(G_{2})$.
\end{enumerate}
\end{proof}

\section{Conclusions and Open Problems}
In this paper, we have defined the independent neighborhood polynomial of a given graph $G$, denoted by $N^{(i)}(G,x)$, as well as the connected neighborhood polynomial of $G$, denoted by $N^{(c)}(G,x)$ which are obtained from the neighborhood polynomial of $G$ after taking particular restrictions into account and then we investigated some of the characteristics of these polynomials.

\vspace{2mm}
\noindent
For a given statement $S$ if $S$ is true $[S]=1$, otherwise $[S]=0$. Let $iso(G)$ be the number of isolated vertices in $G$. We can reformulate the neighborhood polynomial of $G$ as follows: 

\begin{align*}
    N(G,x)&=\sum_{X\in \mathcal{N}(G)}x^{|X|}\\
    &=\sum_{X\in \mathcal{N}(G)}[k(G[X])=|X|]x^{|X|}+\sum_{X\in \mathcal{N}(G)}[k(G[X])=1]x^{|X|}\\
    &+\sum_{X\in \mathcal{N}(G)}[1<k(G[X])<|X|]x^{|X|}\\
    &-1-(n-iso(G))x.
\end{align*}

\noindent
Clearly, $\sum_{X\in \mathcal{N}(G)}[k(G[X])=|X|]x^{|X|}$ is exactly $N^{(i)}(G,x)$, and similarly $\sum_{X\in \mathcal{N}(G)}[k(G[X])=1]x^{|X|}$ is nothing but $N^{(c)}(G,x)$. By $-1-(n-iso(G))x$, we correct the over counting of the empty set and all non-isolated vertices. On the other hand, $\sum_{X\in \mathcal{N}(G)}[1<k(G[X])<|X|]x^{|X|}$, which we denote by $N^{(d)}(G,x)$, generates all vertex subsets of open neighborhoods of vertices in $G$ which induce a disconnected graph with at least one component containing more than one vertex. To conclude, we have
$$N(G,x)=N^{(i)}(G,x)+N^{(c)}(G,x)+N^{(d)}(G,x)-1-(n-iso(G))x.$$

\noindent
Can we find an efficient way or a combinatorial approach to calculate $N^{(d)}(G,x)$ for a graph $G$, given $N^{(i)}(G,x)$ or/and $N^{(c)}(G,x)$?
\section{Acknowledgements}
This work was supported by the \emph{European Social Fund (ESF)}.

\end{document}